\documentclass[11pt,fullpage]{article}
\usepackage{graphicx}
\usepackage{latexsym}
\usepackage{color}

\usepackage{amsmath}
\usepackage{graphicx}
\usepackage{amssymb}
\usepackage{amsthm}
\usepackage{amsfonts}

\renewcommand{\baselinestretch} {1.3}

        \makeatletter 
\def\singlespace{\def\baselinestretch{1}\@normalsize}

\newtheorem{theorem}{Theorem}
\newtheorem{lemma}{Lemma}
\newtheorem*{assumption}{Assumption (P)}
\newtheorem{corollary}{Corollary}
\newtheorem{proposition}{Proposition}

\newcommand{\be}{\begin{equation}}
\newcommand{\ee}{\end{equation}}
\newcommand{\beqn}{\begin{eqnarray}}
\newcommand{\eeqn}{\end{eqnarray}}
\newcommand{\bt}{\begin{theorem}}
\newcommand{\et}{\end{theorem}}
\newcommand{\bl}{\begin{lemma}}
\newcommand{\bp}{\begin{proposition}}
\newcommand{\ep}{\end{proposition}}
\newcommand{\bc}{\begin{corollary}}
\newcommand{\ec}{\end{corollary}}

\newcommand{\fr}[1]{(\ref{#1})} 

\newcommand{\by}{{\bf y}}
\newcommand{\bd}{{\bf d}}
\newcommand{\Jc}{{\cal J}^c_0}
\newcommand{\J}{{\cal J}_0}

\newcommand{\bm}{\mbox{\boldmath $\mu$}}
\newcommand{\btm}{\mbox{\boldmath{$\tilde \mu$}}}
\newcommand{\bhm}{\mbox{\boldmath{$\hat \mu$}}}
\newcommand{\beps}{\mbox{\boldmath $\epsilon$}}

\begin{document}
\title{\Large{\bf Estimation of a sparse group of sparse vectors}}

\author{
{\bf Felix Abramovich}\\
Department of Statistics and Operations Research\\
Tel Aviv University\\
Tel Aviv 69978, Israel\\
felix@post.tau.ac.il
\\ \\
{\bf Vadim Grinshtein}\\
Department of Mathematics\\
The Open University of Israel\\
Raanana 43537, Israel\\
vadimg@openu.ac.il}

\date{}
\maketitle

\vspace{1cm}
\begin{abstract}
We consider a problem of estimating a sparse group of sparse
normal mean vectors. The proposed approach is based on penalized
likelihood estimation with complexity penalties on the number of nonzero
mean vectors and the numbers of their ``significant'' components,
and can be performed by a computationally fast
algorithm. The resulting estimators are developed within Bayesian framework
and can be viewed as MAP estimators.
We establish their adaptive
minimaxity over a wide range of sparse and dense settings.
The presented short simulation study demonstrates the efficiency of
the proposed approach that successfully competes with the recently
developed sparse group lasso estimator.
\end{abstract}


\section{Introduction}
Suppose we observe a series of $m$ independent $n$-dimensional
Gaussian vectors $\by_1,...,\by_m$ with independent components and
common variance: 
\be \label{eq:model} \by_j=\bm_j+\beps_j,\;\;\;
\beps_j  \stackrel{i.i.d.} \sim {\cal N}_n({\bf 0},\sigma_n^2
I_n),\;\;\;j=1,...,m 
\ee 
The variance $\sigma_n^2>0$, which may depend on $n$, is assumed to be
known, and the goal is to estimate the unknown mean vectors
$\bm_1,...,\bm_m$.

The key extra assumption on the model \fr{eq:model} is {\em both} within-
and between-vectors sparsity (hereafter {\em within}- and {\em
between}-sparsity for brevity). More specifically, we assume that
part of $\bm_j$'s are identically zero vectors and the entire information
in the noisy data is contained only in a small fraction of them
(between-sparsity). Moreover, even within nonzero $\bm_j$'s,
most of their components are still zeroes or at least ``negligible''
(within-sparsity). Formally, the within-sparsity can be quantified
in terms of $l_0$, strong or weak $l_p$-balls introduced further.
Neither the indices of non-zero $\bm_j$'s nor the
locations of their ``significant'' components are known in advance.

Such a model appears in the variety of statistical applications as we illustrate
by the following two examples.

\vspace{.5cm} \noindent {\em Example 1. Additive models}. Consider a
nonparametric regression model
$y_i=f(x_{1i},...,x_{m_i})+\epsilon_i,\;i=1,...,n$, where $f:
\mathbb{R}^m \rightarrow \mathbb{R}$ is the unknown regression
function assumed to belong to some class of functions (e.g.,
H\'older, Sobolev or Besov classes), and $\epsilon_i
\stackrel{i.i.d} \sim {\cal N}(0,\sigma_n^2)$. Estimating $f$ in
such a general setup suffers from a severe ``curse of
dimensionality'', where typically the sample size $n$ should grow
exponentially with the dimensionality $m$ to achieve consistent
estimation. It is essential then to place some extra restrictions on
the complexity of $f$. One of the most common approaches is to
consider the additive models, where
$f(x_1,...,x_m)=f_1(x_1)+...+f_m(x_m)$ and each component $f_j$ lies
in some smoothness class. In addition, similar to sparse linear
regression models, it is often reasonable to assume that only part
of predictors among $x_1,...,x_m$ are really ``significant'', while
the impact of others is negligible if at all. Such {\em sparse}
additive models are especially relevant for $m \sim n$ and $m \gg n$
setups and have been considered in Lin \& Zhang (2006), Meier, van
de Geer \& Buhlmann (2009), Ravikumar {\em et al.} (2009), Raskutti,
Wainwright \& Yu (2012).

Expand each $f_j,\;j=1,...,m$ into (univariate) orthonormal series
$\{\psi_{ij}\}$ as $\sum \mu_{ij} \psi_{ij}(x_j)$,
where $\mu_{ij}=\int f_j(x_j) \psi_{ij}(x_j)dx_j$. The original
nonparametric additive model is then transformed into the equivalent
problem of estimating vectors of corresponding coefficients
$\bm_1,...,\bm_m$ within Gaussian noise \fr{eq:model}, where for
sparse additive models, most of $\bm_j$ are zeroes
(between-sparsity). Moreover, for a properly chosen bases
$\{\psi_{ji}\}$ (e.g., Fourier series for Sobolev or wavelets for
more general Besov classes), the nonzero $\bm_j$ will be also sparse
(within-sparsity).

\vspace{.5cm} \noindent {\em Example 2. Time-course microarray
experiments}. In time-course microarray experiments the data
consists of measurements of differences in the expression levels
between ``treated'' and ``control'' samples of $m$ genes recorded at
different times. A record on $j$-th gene at time point $t_i$ is
modelled as a measurement of an (unknown) expression profile function
$f_j(t)$ at time $t_i$ corrupted by Gaussian noise. The expression
of most genes are the same in both groups ($f_j \equiv 0$) and the
goal is to identify the differentially expressed genes and estimate
the corresponding non-identically zero expression profile functions
$f_j$. Similar to the previous example, each $f_j$ is commonly expanded into
some ``parsimonious'' orthonormal basis (e.g., Legendre polynomials,
Fourier or wavelets) as $f_j(t)=\sum_i \mu_{ij} \psi_{ij}(t)$ 
and in the coefficients domain the original functional model becomes
$$
y_{ij}=\mu_{ij}+z_{ij},\;\;\;j=1,...,m;\;i=1,...,n
$$
where $y_{ij}$ are empirical coefficients of the data on $j$-th gene
and $z_{ij}$ are Gaussian noise (see, e.g., Angelini {\em et. al}, 2007).
For most genes, $\bm_j \equiv 0$ (between-sparsity), while due to the
parsimonity of the chosen basis, for differentially expressed genes, $\bm_j$
will still have sparse representation (within-sparsity).

\vspace{.5cm}

To estimate $\bm_1,...,\bm_m$ in \fr{eq:model} under the assumptions of
between- and within-sparsity we proceed as follows. 
From a series of
pioneer works of Donoho \& Johnstone in nineties (e.g., Donoho \& Johnstone,
1994ab),  it is well-known
that the optimal strategy for estimating a {\em single} sparse vector $\bm_j$ from $\by_j$ is thresholding.
Various threshold estimators $\bhm_j$ can be considered
as penalized likelihood estimators, where
$$
\bhm_j=\arg \min_{\btm_j \in \mathbb{R}^n}||\by_j-\btm_j||_2^2+Pen_j(\btm_j),
$$
corresponding to different choices of penalties $Pen_j(\btm)$. In
particular, the $l_1$-type penalty
$Pen_j(\btm_j)=\lambda||\btm_j||_1$ leads to soft thresholding of
components of $\btm_j$ with a constant threshold $\lambda/2$ that
coincides with the lasso estimator of Tibshirani (1996). Wider
classes of penalties on the {\em magnitudes} of components
$\tilde{\mu}_{ij}$ are discussed in Antoniadis \& Fan (2001). In this paper
we consider the
$l_0$ or complexity type penalties $Pen_j(||\btm_j||_0)$ on the
{\em number of nonzero} components $\tilde{\mu}_{ij}$, where
$||\btm_j||_0=\#\{i: \tilde{\mu}_{ij} \ne 0\}$, 
that yield hard thresholding rules. In the simplest case, where
$Pen_j(||\btm_j||_0)=\lambda ||\btm_j||_0$, the resulting (constant)
threshold is $\sqrt{\lambda}$. More general complexity penalties
were studied in Birg\'e \& Massart (2001), Abramovich, Grinshtein \&
Pensky (2007), Abramovich {\em et al.} (2010) and Wu \& Zhou (2012).

Penalizing each $\btm_j$ separately, however, essentially ignores
the between-sparsity, where it is assumed that most of $\bm_j$ are
identically zeroes and should be obviously estimated by $\bhm_j
={\bf 0}$. Thus, simultaneous estimation of all $m$ mean vectors in
\fr{eq:model} should involve an additional penalty $Pen_0(\cdot)$ on the
number of nonzero $\bhm_j$'s that are now defined as solutions of
the following criterion: \be  \label{eq:est} \min_{\btm_1,...,\btm_m
\in \mathbb{R}^n} \left\{\sum_{j=1}^m \left\{
||\by_j-\btm_j||_2^2+Pen_j(||\btm_j||_0)\right\}+ Pen_0(k)\right\},
\ee where $k=\#\{j: \btm_j \ne {\bf 0}\}$. In this paper we
investigate the optimality of such an approach for estimating
$\bm_1,...,\bm_m$ under various within- and between-sparsity setups.
In particular, we specify the classes of complexity penalties
$Pen_j(||\btm_j||_0)$ and $Pen_0(k)$ on respectively within- and
between sparsity for which the resulting estimators
$\bhm_1,...,\bhm_m$ achieve asymptotically minimax rates
simultaneously for the wide range of sparse and dense cases. Such types of
penalties naturally arise within a Bayesian model selection
framework. In this sense, this paper extends the results of Bayesian
MAP testimation approach developed in Abramovich, Grinshtein \&
Pensky (2007) and Abramovich {\em et al.} (2010) for estimating a
single normal mean vector to simultaneous estimation of a group of
$m$ vectors in the model \fr{eq:model}.

It is interesting to compare the proposed complexity penalization \fr{eq:est}
with lasso-type procedures. Similar to $l_0$-type penalization,
the vector-wise use of the original lasso of Tibshirani
(1996) for estimating each $\bm_j$ in \fr{eq:model}
results in per-component (soft) thresholding of each $\by_j$ that
handles within-sparsity but ignores between-sparsity. To address
the latter, Yuan \& Lin (2006) proposed a {\em group lasso} that for the
particular model \fr{eq:model} at hand solves
$$
\min_{\btm_1,...,\btm_m \in \mathbb{R}^n}
\sum_{j=1}^m \left\{||\by_j-\btm_j||_2^2+\lambda ||\btm_j||_2 \right\}
$$
It can be easily shown that in such a setup, the group lasso
estimator is available in the closed form, namely,
$\bhm_j=(1-\frac{\lambda/2}{||{\bf y}_j||_2})_+{\bf
y}_j,\;j=1,...,m$ which is the vector-level ``shrink-or-kill''
thresholding with a threshold $\lambda/2$. The $\bhm_j$'s are,
therefore, either entirely zero or do not have zero components at
all. As a result, the group lasso does not handle within-sparsity. To
combine both types of sparsity, Friedman, Hastie \& Tibshirani
(2010) introduced the {\em sparse group lasso} that for the model
\fr{eq:model} is defined as \be \label{eq:sparselasso}
\min_{\btm_1,...,\btm_m \in \mathbb{R}^n} \sum_{j=1}^m
\left\{||\by_j-\btm_j||_2^2+ \lambda_1 ||\btm_j||_2 + \lambda_2
||\btm_j||_1\right\} \ee yielding
$\bhm_j=(1-\frac{\lambda_1/2}{||\tilde{\bf y}_j||_j})_+\tilde{\bf
y}_j,\;j=1,...,m$, where $\tilde{y}_{ij}={\rm
sign}(y_{ij})(|y_{ij}|-\lambda_2/2)_+,\;i=1,...,n$ is the result of
component-level soft thresholding of each ${\bf y}_j$ with a
threshold $\lambda_2/2$.

To the best of our knowledge, there are no theoretical results on optimality
of sparse group lasso similar to those presented in this paper for the
complexity penalized estimators \fr{eq:est}. Moreover, we believe that, generally,
$l_0$-type penalties are more ``natural'' for representing sparsity and the
main reason for other types of penalties ($l_1$ in particular) are mostly
computational.
For a general regression model, complexity penalties indeed imply
combinatorial search over all possible models, while, for example,
sparse group lasso estimator
can be still efficiently computed by numerical iterative algorithms
(see Friedman, Hastie \& Tibshirani, 2010 and Simon {\em et al.}, 2011 for 
details).
However, for the model
\fr{eq:model}, that can be essentially viewed as a special case of a general
regression setup, \fr{eq:est} can be
also solved by fast algorithms (see Section \ref{sec:bayes}) that makes such computational
arguments irrelevant.

The paper is organized as follows. In Section \ref{sec:bayes} we
develop a Bayesian formalism that gives raise to penalized
estimators \fr{eq:est}. The asymptotic (as both $m$ and $n$ increase)
adaptive minimaxity of the resulting 
{\em sparse group MAP} estimators over various sparse and dense settings is 
investigated in Section \ref{sec:main}. 
The short simulation study is presented in Section
\ref{sec:examples} and some concluding remarks are given in Section
\ref{sec:remarks}. All the proofs are placed in the Appendix.

\section{Bayesian sparse group MAP estimation} \label{sec:bayes}
Consider again the model \fr{eq:model}. If we knew the indices of nonzero
vectors $\bm_j$ and the locations of their
``significant'' entries $\mu_{ij}$, we would evidently estimate them  by 
the corresponding $y_{ij}$ and set others to zero.
Hence, the original problem is
essentially reduced to finding an $n \times m$ indicator matrix $D$, where
$d_{ij}$ indicates whether $\mu_{ij}$ is ``significant'' or not,
and can be viewed as a model
selection problem. Note that due to between- and within-sparsity assumptions,
the matrix $D$ should be sparse in the double sense:
only part of $D$'s columns ${\bf d}_j$ are supposed to be nonzeroes, and
even nonzero columns are sparse.

We introduce first some notations.
Let $\J$ and $\Jc$ be the sets of indices corresponding respectively
to zero and  nonzero mean vectors $\bm_j$'s, and
$m_0=|\Jc|=\#\{j: \bm_j \ne {\bf 0}, \;j=1,...,m\}$.
Denote by $h_j=\sum_{i=1}^n d_{ij}=
\#\{i: \mu_{ij} \ne 0,\;i=1,...,n\}$ the number of nonzero components in
$\bm_j$, where evidently $h_j=0$ for $j \in \J$.

Consider the following Bayesian model selection procedure for
identifying nonzero components $\mu_{ij}$ or, equivalently, the
indicator matrix $D$. To capture the between- and within-sparsity
assumptions we place a hierarchical prior on $D$.
We first assume some prior distribution on the number of nonzero
mean vectors $m_0 \sim \pi_0(m_0)>0,\;m_0=0,...,m$. For a given $m_0$,
assume that all ${m \choose m_0}$ different configurations of zero
and nonzero mean vectors are equally likely, that is, conditionally
on $m_0$,
$$
P(\Jc\; \bigl|\;|\Jc|=m_0)={m \choose m_0}^{-1}
$$
Obviously, $h_j\bigl|\{j \in \J\} \sim \delta(0)$ and, thus, ${\bf
d}_j\bigl|\{j \in \J\} \sim \delta({\bf 0})$ and $\bm_j\bigl|\{j \in
\J\} \sim \delta({\bf 0})$. For nonzero $\bm_j$ we place independent
priors $\pi_j(\cdot)$ on the number of their nonzero components,
that is, $h_j\bigl|\{j \in \Jc\} \sim \pi_j(h_j)>0,\;h_j=1,...,n$. In
this case, we again assume that for a given $h_j$, all
possible ${n \choose h_j}$ indicator vectors ${\bf d}_j$ with
$h_j$ nonzero components
have the same prior probabilities and, therefore,
$$
P({\bf d}_j \; \bigl|\; ||{\bf d}_j||_0=h_j,j \in \Jc)={n \choose h_j}^{-1}
$$
Finally, to complete the prior for \fr{eq:model}, we have
$\mu_{ij}\bigl|d_{ij}=0 \sim \delta(0)$, while
nonzero $\mu_{ij}$ are assumed to be i.i.d. $N(0,\gamma \sigma_n^2)$,
where $\gamma>0$.

A straightforward Bayesian calculus yields the posterior probability
for a given indicator matrix $D$:
$$
P(D\bigl|\by) \propto \pi_0(m_0) {m \choose m_0}^{-1} \prod_{j \in
\Jc} \left\{\pi_j(h_j){n \choose h_j}^{-1}
(1+\gamma)^{-\frac{h_j}{2}} e^{\frac{\gamma}{\gamma+1}
\frac{\sum_{i=1}^n y_{ij}^2 d_{ij}}{2\sigma_n^2}} \right\}
$$
Given the
posterior distribution $P(D|\by)$ we apply the maximum {\em a
posteriori} (MAP) rule to choose the most likely configuration of
zero and nonzero $\mu_{ij}$ that leads to the following MAP
criterion: 
\be 
\sum_{j \in \Jc} \left\{\sum_{i=1}^n
y_{ij}^2 d_{ij}+ 2\sigma_n^2(1+1/\gamma) \ln\left(\pi_j(h_j){n
\choose h_j}^{-1}(1+\gamma)^{-\frac{h_j}{2}}\right)\right\}+
2\sigma_n^2(1+1/\gamma)\ln\left(\pi_0(m_0){m \choose
m_0}^{-1}\right) \rightarrow \max_D 
\label{eq:map0} 
\ee 
From \fr{eq:map0} it follows
immediately that for a given $h_j>0$ the optimal choice $\hat{\bf
d}_j(h_j)$ for ${\bf d}_j$ is $\hat{d}_{ij}(h_j)=1$ for the $h_j$
largest $|y_{ij}|$ and zero otherwise. The criterion \fr{eq:map0} is
then reduced to \be \label{eq:map1} \sum_{j \in \Jc}
\left\{\sum_{i=1}^{h_j} y_{(i)j}^2+ 2\sigma_n^2(1+1/\gamma)
\ln\left(\pi_j(h_j){n \choose
h_j}^{-1}(1+\gamma)^{-\frac{h_j}{2}}\right)\right\} +
2\sigma_n^2(1+1/\gamma)\ln\left(\pi_0(m_0){m \choose
m_0}^{-1}\right) \rightarrow \max_D, \ee where $|y_{(1)j}| \geq ...
\geq |y_{(n)j}|$. For every $j=1,...,m$ define
\begin{eqnarray}
\hat{h}_j& = &\arg \min_{1 \leq h_j \leq n} \left\{\sum_{i=h_j+1}^n
y_{(i)j}^2+2\sigma_n^2(1+1/\gamma) \ln\left(\pi^{-1}_j(h_j){n
\choose h_j}(1+\gamma)^{\frac{h_j}{2}}\right)\right\}
\nonumber \\
& = & \arg \min_{1 \leq h_j \leq n}
\left\{-\sum_{i=1}^{h_j}y_{(i)j}^2+2\sigma_n^2(1+1/\gamma)\ln\left(\pi^{-1}_j(h_j){n
\choose h_j}(1+\gamma)^{\frac{h_j}{2}}\right)\right\} \label{eq:hj}
\end{eqnarray}
Then, \fr{eq:map1} is equivalent to minimizing \be \sum_{j \in \Jc}
\left\{-\sum_{i=1}^{\hat{h}_j} y^2_{(i)j}+2\sigma_n^2(1+1/\gamma)
\ln\left(\pi^{-1}_j(\hat{h}_j){n \choose \hat{h}_j}
(1+\gamma)^{\frac{\hat{h}_j}{2}}\right)\right\}+
2\sigma_n^2(1+1/\gamma)\ln\left(\pi^{-1}_0(m_0){m \choose
m_0}\right) \label{eq:map2} \ee over all subsets of indices $\J
\subseteq \{1,...,m\}$. Define \be \label{eq:wj}
W_{j}=-\sum_{i=1}^{\hat{h}_j} y^2_{(i)j}+2\sigma_n^2(1+1/\gamma)
\ln\left(\pi^{-1}_j(\hat{h}_j){n \choose
\hat{h}_j}(1+\gamma)^{\frac{\hat{h}_j}{2}}\right) \ee Then,
\fr{eq:map2} is obviously reduced to \be \label{eq:map3} \min_{0
\leq m_0 \leq m}\left\{ \sum_{j=1}^{m_0}W_{(j)} +
2\sigma_n^2(1+1/\gamma)\ln\left(\pi^{-1}_0(m_0){m \choose
m_0}\right) \right\}, \ee where $W_{(1)} \leq ... \leq W_{(m)}$ and
for $m_0=0$ the sum in the RHS of \fr{eq:map2} evidently does not
appear.

Summarizing, the efficient simple algorithm for finding the proposed
sparse group MAP estimators of $\bm_1,...,\bm_m$ in \fr{eq:model}
can be formulated as follows:

\vspace{.5cm}
{\em
\centerline{\bf Sparse group MAP estimation algorithm}

\begin{enumerate}
\item For every $j=1,...,m$, find $\hat{h}_j$ in \fr{eq:hj} and calculate
the corresponding $W_j$ in \fr{eq:wj}.

\item Order $W_j$ in ascending order $W_{(1)} \leq ... \leq W_{(m)}$ and
find
$$
\hat{m}_0=\arg\min_{0 \leq m_0 \leq m}\left\{
\sum_{j=1}^{m_0}W_{(j)} +
2\sigma_n^2(1+1/\gamma)\ln\left(\pi^{-1}_0(m_0){m \choose
m_0}\right) \right\}
$$

\item Let $\hat{\Jc}$ be the set of indices corresponding to the $\hat{m}_0$
smallest $W_j$. Set $\bhm_j \equiv {\bf 0}$ for all $j \in \hat{\J}$, while
for $j \in \hat{\Jc}$, take the $\hat{h}_j$ largest $|y_{ij}|$ and threshold
others, that is, $\hat{\mu}_{ij}=y_{ij} \mathbb{I}\{|y_{ij}| \geq |y_{(\hat{h}_j)j}|\},\;i=1,...,n,\;j \in \hat{\Jc}$, where $|y_{(1)j}| \geq ... \geq |y_{(n)j}|$.
\end{enumerate}
}
\vspace{.5cm}

The resulting estimation procedure combines therefore vector-wise
and component-wise thresholding. It is easily verified that the
minimizer of (\ref{eq:map2}) is, in fact, the penalized likelihood
estimator \fr{eq:est} with the complexity penalties \be
\label{eq:penaltyj} Pen_j(0)=0,\;Pen_j(h_j)= 2\sigma_n^2(1+1/\gamma)
\ln\left(\pi^{-1}_j(h_j){n \choose
h_j}(1+\gamma)^{\frac{h_j}{2}}\right), \;h_j=1,...,m \ee and \be
\label{eq:penalty0}
Pen_0(m_0)=2\sigma_n^2(1+1/\gamma)\ln\left(\pi^{-1}_0(m_0){m \choose
m_0}\right),\;m_0=0,...,m
\ee

The specific types of penalties $Pen_j(\cdot)$'s and $Pen_0(\cdot)$
depend on the choices of priors $\pi_j(\cdot)$'s and $\pi_0(\cdot)$.
For example, binomial priors $m_0 \sim B(m,\xi_0)$ and $h_j \sim
B(n,\xi_j)$ yield linear type penalties $Pen(m_0)=2\sigma_n^2
\lambda^2_0 m_0$ and  $Pen_j(h_j)=2\sigma_n^2 \lambda^2_j h_j$
respectively, where $\lambda^2_0=(1+1/\gamma)\ln\{(1-\xi_0)/\xi_0\}$
and $\lambda^2_j=(1+1/\gamma)\ln\{\sqrt{1+\gamma}(1-\xi_j)/\xi_j\}$.
For such a choice of $\pi_j(\cdot)$, $W_j$ in \fr{eq:wj} is
essentially obtained by hard thresholding of $\by_j$ with a {\em
constant} threshold $\sqrt{2}\sigma_n\lambda_j$.  In particular,
$\xi_j=\sqrt{\gamma+1}/(\sqrt{\gamma+1}+n^{\gamma/(\gamma+1)})$
leads to the universal thresholding of Donoho \& Johnstone (1994a)
with $\lambda_j=\sqrt{\ln n}$. The (truncated) geometric priors
$\pi_j(h_j) \propto q_j^{h_j},\;h_j=1,...,n$ for some $0 < q_j < 1$,
imply the (nonlinear) so-called $2k\ln(n/k)$-type penalties. The
optimality of the resulting hard thresholding estimator with a {\em
data-driven} threshold for estimating a single normal mean vector
has been shown in Abramovich, Grinshtein \& Pensky (2007),
Abramovich {\em et. al} (2010), Wu \& Zhou (2012).

\section{Adaptive minimaxity of sparse group MAP estimators} \label{sec:main}
In this section we investigate the goodness of the proposed sparse
group MAP estimators \fr{eq:est} with the penalties
\fr{eq:penaltyj}-\fr{eq:penalty0}, where 
the goodness-of-fit is measured by the global quadratic risk
$\sum_{j=1}^m E||\bhm_j-\bm_j||^2_2$. 
We establish their
asymptotic minimaxity over a wide range of sparse and dense settings.
To derive these results we need the following assumption on the priors
$\pi_j(\cdot)$:
\begin{assumption}
\label{as:P}
Assume that
\begin{equation}
\pi_j(h) \leq {n \choose h}e^{-c(\gamma)h},\;h=1,...,n,\;j=1,...,m, \label{eq:P}
\end{equation}
where $c(\gamma)=8(\gamma+3/4)^2>9/2$.
\end{assumption}
Assumption (P) is, in fact, not restrictive. Indeed,
the obvious inequality ${n \choose h} \geq (n/h)^h$ implies that
for {\em any} $\pi_j(\cdot)$, \fr{eq:P} holds
for all $h \leq ne^{-c(\gamma)}$. In particular, Assumption (P) is satisfied for
binomial priors $B(n,\xi_j)$ with $\xi_j \leq e^{-c(\gamma)}/(1+e^{-c(\gamma)})$
and (truncated) geometric priors.

First, we obtain a general upper bound for the quadratic risk
of the sparse group MAP estimator that will be the key
for deriving its asymptotic minimaxity.

\begin{theorem}[general upper bound] \label{th:upper}
Consider the sparse group
MAP estimators $\bhm_1,....,\bhm_m$ \fr{eq:est}
of $\bm_1,...,\bm_m$ with the complexity penalties \fr{eq:penaltyj}-\fr{eq:penalty0}
in the model \fr{eq:model}.
Under Assumption (P) we have
\begin{eqnarray}
\sum_{j=1}^m E||\bhm_j-\bm_j||^2_2 & \leq & c_1(\gamma) \min_{\J \subseteq \{1,...,m\}} \left\{
\sum_{j \in \Jc}
\min_{1 \leq h_j \leq n}
\left(\sum_{i=h_j+1}^n \mu^2_{(i)j}+Pen_j(h_j)\right) \right. \nonumber \\
& + & \left. \sum_{j \in \J}\sum_{i=1}^n
\mu^2_{ij}+Pen_0(|\Jc|)\right\} +c_2(\gamma) \sigma_n^2 (1-\pi_0(0)),
\label{eq:upper}
\end{eqnarray}
where $|\mu_{(1)j}| \geq ... \geq |\mu_{(n)j}|$ and $c_1(\gamma)$, $c_2(\gamma)$
depend only on $\gamma$.
\end{theorem}

The results of Theorem \ref{th:upper} hold for any
normal mean vectors $\bm_1,...,\bm_m$.
Now we consider \fr{eq:model} under the extra within- and between-sparsity assumptions
that will be defined more rigorously below.

The between-sparsity is naturally measured by the number $m_0$ of
nonzero $\bm_j$'s. The within-sparsity can be introduced in several
ways. The most intuitive measure of within-sparsity of a single
normal mean vector $\bm \in \mathbb{R}^n$ is the number of its
nonzero components, that is, its $l_0$ quasi-norm $||\bm||_0$.
Define then an $l_0$-ball $l_0[\eta]$ of standardized radius $\eta$
as a set of $\bm$ with at most a proportion $\eta$ of non-zero
entries, that is
$$
l_0[\eta]=\{\bm \in \mathbb{R}^n~: ||\bm||_0 \leq \eta n \}
$$
One can argue that in many practical settings, it is more reasonable
to assume that the components $\mu_i$'s of $\bm$ are not exactly
zero but ``small''. In a wider sense the within-sparsity of $\bm$
can be then defined by the proportion of its large entries.
Formally, define a weak $l_p$-ball $m_p[\eta]$ with a standardized
radius $\eta$ as
$$
m_p[\eta]=\{\bm \in \mathbb{Re}^n~: |\mu|_{(i)} \leq \sigma_n \eta
(n/i)^{1/p},\;i=1,...,n\},
$$
where $\mu_{(1)} \geq ... \geq \mu_{(n)}$ are the ordered components of $\bm$.
For $\bm \in m_p[\eta]$, the proportion of $|\mu_i|$'s larger 
than $\sigma_n \delta$ for some $\delta>0$ is at
most $(\eta/\delta)^p$.

Within-sparsity can be also measured in terms of the $l_p$-norm of $\bm$, where
a strong $l_p$-ball $l_p[\eta]$ with standardized radius $\eta$ is defined as
$$
l_p[\eta]=\{\bm \in \mathbb{Re}^n~: \frac{1}{n}\sum_{i=1}^n|\mu_i|^p
\leq \sigma_n^p \eta^p\}
$$
There are well-known relationships between these types of balls.
The $l_p$-norm approaches $l_0$ as $p$ decreases, while a weak $l_p$-ball
contains the corresponding strong $l_p$-ball but only just:
$$
l_p[\eta] \subset m_p[\eta] \not\subset l_{p'}[\eta],\;p'>p
$$

We recall first the known results on minimax rates for estimating a
{\em single} normal mean vector $\bm$ over different types of balls
introduced above. Let $\Theta[\eta_n] \subset \mathbb{R}^n$ be any
of $l_0[\eta_n],l_p[\eta_n]$ or $m_p[\eta_n]$, where the
standardized radius $\eta$ might depend on $n$. The corresponding
minimax quadratic risk for estimating a single $\bm$ ($m=1$) over
$\Theta[\eta_n]$ in \fr{eq:model} is
$R(\Theta[\eta_n])=\inf_{\btm}\sup_{\bm \in
\Theta[\eta_n]}E||\btm-\bm||^2_2$, where the infimum is taken over
all estimates $\btm$ of $\bm$. For $p>0$ define
$\eta_{0n}=n^{-1/\min(p,2)} \sqrt{\ln n}$. Depending
on the behaviour of $\eta_n$ as $n$ increases, we distinguish between
three cases for $p>0$ and two cases for $p=0$: 
\begin{itemize}
\item[a)] {\em dense}, where
$\eta_n \not \rightarrow 0$
\item[b)] {\em sparse}, where $\eta_n \rightarrow 0$ but $\eta_n/\eta_{0n} \not \rightarrow 0$ for $p>0$ and, obviously, $\eta_n \geq n^{-1}$ for $p=0$
\item[c)] {\em super-sparse} (for $p>0$), where $\eta_n/\eta_{0n} \rightarrow 0$\end{itemize} 
The corresponding minimax convergence rates over $R(\Theta[\eta_n])$ for
various cases and $p$ are summarized in Table \ref{tab:rates} below
(see Donoho {\em et. al}, 1992; Johnstone, 1994; Donoho \&
Johnstone, 1994b). 

The rates for $m_p[\eta_n]$ are the same as for $l_p[\eta_n]$
except $p=2$, where there is an additional log-term.
Table \ref{tab:rates} defines dense and sparse
zones for $p=0$ and $p \geq 2$, and dense, sparse and super-sparse
zones for $0 < p < 2$ of different minimax rates.

\begin{table}[htb]
\begin{center}
\begin{tabular}{|l|c|c|c|}
\hline
Case  & $p=0$ & $0 < p < 2$ & $p \geq 2$ \\
\hline

dense case & $\sigma_n^2 n$ & $\sigma_n^2 n$ &
 $\sigma_n^2 n$ \\

sparse case & $\sigma_n^2 n \eta_n(\ln \eta_n^{-1})$ & 
$\sigma_n^2 n \eta_n^p(\ln \eta_n^{-p})^{1-p/2}$ &
$\sigma_n^2 n \eta_n^2$ \\

super-sparse case &
$ - $  &
$\sigma_n^2n^{2/p}\eta_n^2$ & $\sigma_n^2 n \eta_n^2$ \\
\hline
\end{tabular}
\caption{Minimax rates (up to multiplying constants) over various 
$l_0[\eta_n]$, $l_p[\eta_n]$ and
$m_p[\eta_n]$-balls. The rates are the same for $l_p[\eta_n]$ and $m_p[\eta_n]$
except $p=2$, where for $m_p[\eta_n]$ there appears the additional log-term
which is not presented in Table \ref{tab:rates} for brevity.}
\label{tab:rates}
\end{center}
\end{table}

Consider now the model \fr{eq:model} for $m \geq 1$. Recall that
$m_0=\#\{j: \bm_j \not \ne {\bf 0}\}$ and $\Jc$ is the set of
indices for nonzero $\bm_j$. In what follows we assume that
$\bm_j \in \Theta_j[\eta_{jn}]$ for $j \in \Jc$, where the types ($l_0$,
weak $m_p$ or strong $l_p$) and the parameters $p$ of the
corresponding balls are not necessarily the same for all $j$.
Furthermore, we allow the priors $\pi_0(\cdot)$ and $\pi_j(\cdot)$
to depend respectively on $m$ and $n$.

Theorem \ref{th:upperth} below defines the asymptotic upper bounds
for the quadratic risks of the sparse group MAP estimator in
\fr{eq:model} under within- and between sparsity assumptions:
\begin{theorem}[upper bounds over sparse and dense settings] \label{th:upperth}
Consider the model \fr{eq:model}, where $\Jc \neq \emptyset$ (not
pure noise). Assume that $\bm_j \in \Theta_j[\eta_{jn}]$ for all $j \in \Jc$,
where $\eta_{jn} \geq n^{-1/\min(p_j,2)} \sqrt{\ln n}$ for all $p_j>0$ 
(excluding, thus, super-sparse cases).

Let $\bhm_1,...,\bhm_m$
be the sparse group MAP estimators \fr{eq:est} with the complexity penalties
\fr{eq:penaltyj}-\fr{eq:penalty0}, where assume that
there exist constants $c_0, c_1>0$ and $c_2>c(\gamma)$ such that
\begin{enumerate}
\item $\pi_0(k) \geq (k/m)^{c_0 k},\;k=1,...,\lfloor m/e\rfloor$ and
$\pi_0(m) \geq e^{-c_0 m}$
\item for all $j=1,...,m$, $\pi_j(\cdot)$ satisfy Assumption (P) and, in
addition,
$\pi_j(h) \geq (h/n)^{c_1 h},\;h=1,...,\lfloor ne^{-c(\gamma)}\rfloor$;
$\;\;\pi_j(n) \geq e^{-c_2 n}$
\end{enumerate}

Then, for any $\Jc \subseteq \{1,...,m\}$ with $|\Jc|=m_0$ and all
$\Theta_j[\eta_{jn}],\; j \in \Jc$, 
\be \label{eq:rate} 
\sup_{\bm_j
\in \Theta_j[\eta_{jn}],j \in \Jc}\sum_{j=1}^m E||\bhm_j-\bm_j||^2_2
\leq C_1(\gamma) \max\left(\sum_{j \in \Jc}
R(\Theta_j[\eta_{jn}]),\sigma_n^2 m_0 \ln(m/m_0)\right) 
\ee 
for some constant $C_1(\gamma)$ depending only on $\gamma$, where the
corresponding $R(\Theta_j[\eta_n])$ are given in Table
\ref{tab:rates} (up to multiplying constants).
\end{theorem}

Theorem \ref{th:upperth} shows that as both $m$ and $n$ increase,
the asymptotic convergence rates in \fr{eq:rate} are either of order
$\sum_{j \in \Jc}R(\Theta_j[\eta_{jn}])$ or $\sigma_n^2 m_0
\ln(m/m_0)$. The former is associated with the optimal rates of
estimating $m_0$ single sparse vectors in $\Theta_j[\eta_{jn}],\; j
\in \Jc$, while the latter appears in the optimal rates in the model
selection and corresponds to the error of selecting a subset of
$m_0$ nonzero elements out of $m$ (see, e.g. Abramovich \&
Grinshtein, 2010; Raskutti, Wainwright \& Yu, 2011; Rigollet \&
Tsybakov, 2011). From Table \ref{tab:rates} it follows that for all
within-dense and within-sparse cases, $C_1 \sigma^2_n \ln n \leq R(\Theta_j [\eta_{jn}]) \leq C_2 \sigma^2_n n,\;j \in \Jc$ for some $C_1,\;C_2 > 0$
and, therefore,
the first term $\sum_{j \in {\cal J}_0^c}R(\Theta_j[\eta_n])$ in the upper bound (\ref{eq:rate})
is always dominating for $m_0 > m/n$, while the second term 
$\sigma^2_n m_0\ln(m/m_0)$ is necessarily the main one for $m_0 < m/e^n$.

One can easily verify that the conditions on the priors 
$\pi_0(\cdot)$ and $\pi_j(\cdot)$ required in Theorem \ref{th:upperth}
are satisfied, for example, for
the (truncated) geometric priors (see Section \ref{sec:bayes}). On
the other hand, no binomial priors $\pi_0=B(m,\xi_0)$ or
$\pi_j=B(n,\xi_j)$ can satisfy all of them: the requirement
$\pi_j(n)=\xi_j^n \geq e^{-c_2n}$ yields $\xi_j \geq e^{-c_2}$,
while to have $\pi_j(1)=n\xi_j(1-\xi_j)^{n-1} \geq n^{-c_1}$ one
needs $\xi_j \rightarrow 0$ as $n$ increases.

\noindent
\newline
\vspace{.1cm} To establish the corresponding lower bound for the
minimax risk, for simplicity of exposition we consider only the two
cases, where $p_j$ for $j \in \Jc$ are either all zeroes or all
positive. In fact, these are the two main scenarios appearing in
various setups.
Somewhat similar results for minimax lower bounds in the particular
context of sparse nonparametric additive models (see Introduction)
appear in Raskutti, Wainwright and Yu (2012).

\begin{theorem}[minimax lower bounds for $l_0$-balls] \label{th:lower1}
Consider the model \fr{eq:model}, where $\bm_j \in l_0[\eta_{jn}],
\; j \in \Jc$. Assume that $|\Jc|=m_0>0$. Then, there exists a
constant $C_2>0$ such that 
\be \label{eq:lower1}
\inf_{\btm_1,...,\btm_m} \sup_{\bm_j \in l_0[\eta_{jn}],j \in
\Jc}\sum_{j=1}^m E||\btm_j-\bm_j||^2_2 \geq C_2 \max\left(\sum_{j
\in \Jc} R(l_0[\eta_{jn}]),\sigma_n^2 m_0 \ln(m/m_0)\right), 
\ee
where the infimum is taken over all estimates $\btm_1,...,\btm_m$ of
$\bm_1,...,\bm_m$.
\end{theorem}
Theorem \ref{th:lower1} shows that, as $m$ and $n$ increase, the
rates in \fr{eq:rate} cannot be improved for $l_0$-balls. 
The proposed sparse group
MAP estimator in this case is, therefore,  adaptive to the unknown degrees of
within- and between-sparsity and is simultaneously rate-optimal (in
the minimax sense) over entire range of dense and sparse $l_0$-balls settings.

The analysis of the case $p_j>0$ is slightly more delicate. 
Note first that due to the embedding properties of $l_p$-balls
for $p>0$ (see above), it is sufficient to establish the minimax lower
bounds for strong $l_p$-balls settings.

\begin{theorem}[minimax lower bounds for $l_p$-balls] \label{th:lower2}
Consider the model \fr{eq:model}, where $\bm_j \in
l_{p_j}[\eta_{jn}], \; j \in \Jc$ and $|\Jc|=m_0>0$.  In addition,
assume that $\eta^2_{jn} \geq n^{-2/\min(p_j,2)} \max\left(\ln n,\ln(m/m_0)\right)$.
Under this additional constraint, there exists a constant $C_2>0$ such that
\be 
\label{eq:lower2} 
\inf_{\btm_1,...,\btm_m} \sup_{\bm_j \in
l_{p_j}[\eta_{jn}],j \in \Jc}\sum_{j=1}^m E||\btm_j-\bm_j||^2_2 \geq
C_2 \max\left(\sum_{j \in \Jc} R(l_{p_j}[\eta_{jn}]),\sigma_n^2 m_0
\ln(m/m_0)\right), \ee 
where the infimum is taken over all
estimates $\btm_1,...,\btm_m$ of $\bm_1,...,\bm_m$.
\end{theorem}

Similar to Theorem \ref{th:lower1}, Theorem \ref{th:lower2} implies
simultaneous optimality (in the minimax sense) of MAP sparse group estimator
over strong and weak $l_p$-balls but with the restriction on $\eta_{jn}$ and $m_0$. In particular, it does not
cover settings with within-super-sparsity but might also exclude part of the 
corresponding within-sparse zone (depending on $m_0$). 
Within- and between-sparsity cannot be ``too strong'' {\em both}.
In fact, the condition  
$\eta^2_{jn} < n^{-2/\min(p_j,2)} \max\left(\ln n,\ln(m/m_0)\right),\;
j \in \Jc$ can be viewed as an extended definition of super-sparsity
for $m>1$.
For such a super-sparse case, the minimax bound (\ref{eq:lower2}) 
does not hold and can be reduced.
Indeed,  
consider the trivial zero estimator $\btm \equiv {\bf 0},\;j=1,...,m$, 
where, evidently, 
\be \label{eq:zero}
\sup_{\bm_j \in 
l_{p_j}[\eta_{jn}],j \in \Jc}\sum_{j=1}^m E||\btm_j-\bm_j||^2_2 
=\sup_{\bm_j \in 
l_{p_j}[\eta_{jn}],j \in \Jc}\sum_{j \in \Jc} ||\bm_j||^2_2 
\ee
The least favourable sequences that maximize $||\bm_j||^2_2$ over 
$l_{p_j}[\eta_{jn}]$ are $(\sigma_n\eta_{jn},...,\sigma_n\eta_{jn})'$
and $(\sigma_n\eta_{jn}n^{1/p_j},0,...,0)'$ for $p_j \geq 2$ and $0<p_j<2$
respectively. Thus, 
$\sup_{\bm_j \in l_{p_j}[\eta_{jn}]}||\bm_j||^2_2=\sigma_n^2\eta_{jn}^2n^{2/\min(p_j,2)}$ 
and the RHS of \fr{eq:zero} is less than $\sigma^2_n m_0 \ln(m/m_0)$  
for $\eta^2_{jn}<n^{-2/\min(p_j,2)} \ln(m/m_0),\;j \in \Jc$. 
This goes along the lines with the corresponding results for estimating
a single normal mean vector, where a zero estimator is known to be 
rate-optimal for the super-sparse case (Donoho \& Johnstone, 1994b).

\section{Simulation study} \label{sec:examples}
A short simulation study was carried out to demonstrate the
performance of the proposed approach.

The data was generated according to the model \fr{eq:model} with
$m=10$ vectors $\bm_j$'s of length $n=100$. Five $\bm_j$'s were
identically zeroes, while the other five had respectively $100, 70,
50, 20$ and $5$ nonzero components randomly sampled from
$N(0,\tau^2),\;\tau=1,3,5$ and zero others. Such a setup covers various
types of within-sparsity. Finally, the
independent standard Gaussian noise $N(0,1)$ was added to all
components of each $\bm_j$. 

We tried binomial and truncated geometric priors for sparse group
MAP estimators. 
For the binomial prior, we performed component-wise
universal hard thresholding of Donoho \& Johnstone (1994a) with a threshold
$\lambda=\sigma\sqrt{2\log n}$ within each vector that essentially corresponds 
to $\xi_j=\sqrt{\gamma+1}/(\sqrt{\gamma+1}+n^{\gamma/(\gamma+1)})$,
where $\gamma=\tau^2/\sigma^2$ (see Section \ref{sec:bayes}),
and used $\xi_0=1/m$. For the geometric prior
we set $q_0=q_j=0.3$. In addition, we compared the
performances of sparse group MAP estimators with the sparse group
lasso estimator \fr{eq:sparselasso} of Friedman, Hastie \&
Tibshirani (2010) described in Introduction. They do not discuss the
optimal choices for $\lambda_1$ and $\lambda_2$ in
\fr{eq:sparselasso}. Some heuristical arguments are given in Simon {\em et al.} (2011). 
In our simulation study we considered instead two oracle-based choices for these
tuning parameter giving thus a significant handicap to sparse group lasso estimators.
Since in simulation examples the true mean vectors $\bm_j$ are known, they
can be used for optimal choosing $\lambda_1$ and $\lambda_2$.
In particular, we considered  a ``semi-oracle'' sparse group lasso
estimator, where we set $\lambda_2=2\sigma \sqrt{2\log n}$ yielding
universal soft thresholding within each vector (see
Introduction) to compare the sparse group lasso with the binomial
sparse group MAP. $\lambda_1$ was chosen by minimizing the 
mean squared error $\sum_{j=1}^m E||\hat{\bm}_j(\lambda_1)- \bm_j||^2_2$ 
estimated by averaging over a series of 1000 replications for each value
of $\lambda_1$ by a grid search. 
In addition, we applied a ``fully oracle'' sparse group lasso estimator, where
both $\lambda_1$ and $\lambda_2$ were chosen to minimize the mean squared error
by the two-dimensional grid. It can be considered as a benchmark for the
performance of sparse group lasso. Table \ref{tab:lambda} provides the resulting
oracle choices for $\lambda_1$ and $\lambda_2$.

\begin{table}[htb]
\begin{center}
\begin{tabular}{|c|c|c|}
\hline
$\gamma$ & $\lambda_1$ & $\lambda_2$ \\
\hline
1        &   11.8      &     0.9     \\
9        &   7.2       &     1.1     \\
25       &   4.7       &     1.3     \\
\hline
\end{tabular}
\caption{The oracle choices for the parameters of the fully oracle
sparse group lasso estimator ($\gamma=\tau^2/\sigma^2$).}
\label{tab:lambda}
\end{center}
\end{table}
Table \ref{tab:lambda} shows that for all $\gamma$, the oracle 
choice for $\lambda_2$ in the sparse group lasso 
is much less than the conservative universal threshold 
$2\sigma \sqrt{2\log n} \approx 6.06$. The oracle thresholding within each 
vector is thus much less severe and keeps more coefficients. 
The oracle choices for $\lambda_1$ were also quite small and, as a result,
for any $\gamma$, no single vector
was thresholded by a fully oracle sparse group lasso, that is, all $\hat{\bm}_j \ne 0$. 
Thus it was really a non-sparse estimator for the considered setup.  

In Table \ref{tab:mse} we present the mean squared errors averaged
over 1000 replications for the four sparse group estimators with the
corresponding standard errors for various $\gamma$ (or, equivalently,
$\tau$). 

\begin{table}[htb]
\begin{center}
\begin{tabular}{|c|c|c|c|c|}
\hline
$\gamma$ & Sparse Group MAP & Sparse Group MAP & Sparse Group Lasso & Sparse Group Lasso \\
       &   (binomial)     &   (geometric)    & (semi-oracle)          & (fully oracle)         \\
\hline
1      &    247.40        &     245.46       &   236.85  &   161.89  \\
       &    (0.71)        &     (0.70)       &   (0.65)  &   (0.43)  \\
\hline
9      &    608.02        &     378.87       &   1120.99 &   403.76   \\
       &    (1.96)        &     (1.20)       &   (2.29)  &    (0.91)  \\
\hline
25      &    549.77        &     351.52       &   1595.91 &   475.47  \\
       &    (1.68)        &     (1.30)       &   (2.79)  &    (1.07)  \\
\hline
\end{tabular}
\caption{MSEs averaged over 1000 replications for four sparse group
estimators and the corresponding standard errors (in brackets) for various $\gamma$.}
\label{tab:mse}
\end{center}
\end{table}

For small $\gamma$ only few largest nonzero components can be
distinguished from the noise that essentially corresponds to a
sparse setting and explains good performance of binomial sparse
group MAP and semi-oracle sparse group lasso estimators based on universal
(respectively, hard and soft) thresholding within each vector. 
For larger $\gamma$, it becomes ``over-conservative''. 
The negative effect of its conservativeness is 
much stronger for the soft than for hard thresholding (see comments below).
The fully oracle sparse group lasso estimator strongly outperforms its 
semi-oracle counterpart especially for $\gamma=9\;(\tau=3)$ and 
$\gamma=25\;(\tau=5)$ also indicating
that the universal thresholding is far from being optimal for sparse group lasso especially 
for moderate and large $\gamma$ (see also our previous comments on the
optimal choice of $\lambda_2$). 

On the other hand, geometric sparse group MAP estimator corresponding
to a nonlinear $2k\ln(n/k)$-type penalty (see Section \ref{sec:bayes})
provides good results for all $\gamma$ nicely following the theoretical results
of Section \ref{sec:main}. Moreover, for $\gamma=9$ and $\gamma=25$,
it outperforms even the fully oracle sparse group lasso estimator that
was essentially thought as a benchmark rather than a fair competitor. 
This indicates that that sparse group lasso faces general problems.
In fact, it may be not so
surprising since soft ``shrink-or-kill'' thresholding inherent for sparse group
lasso is well-known to be superior to hard ``keep-or-kill'' thresholding
in sparse group MAP estimation for 
small coefficients but worse for large ones due to the additional
shrinkage. Moreover, sparse group lasso essentially involves a {\em double}
amount of
shrinkage - both within vectors and at each entire vector as a whole (see
\fr{eq:sparselasso}).
It thus causes
unnecessary extra bias growing with $\gamma$ that outweighs the benefits of variance reduction.
Similar phenomenon appears also for na\"ive elastic set estimation (Zou \& Hastie,
2005).

\section{Concluding remarks} \label{sec:remarks}
In this paper we considered estimation of a sparse group of sparse
normal mean vectors. The proposed approach is based on penalized
likelihood estimation with complexity penalties on both between- and
within-sparsity and can be performed by a computationally fast
algorithm. The resulting estimators naturally arise within Bayesian
framework and can be viewed as MAP estimators corresponding to the
priors on the number of nonzero mean vectors and the numbers of their
nonzero components. Such a Bayesian perspective provides a natural tool
for obtaining a wide class of penalized likelihood estimators with various
complexity penalties.

We established the adaptive minimaxity of sparse group MAP estimators
to the unknown degree of between- and
within-sparsity over a wide range of sparse and dense settings.
The short simulation study demonstrates the efficiency of
the proposed approach that outperforms the recently
presented sparse group lasso estimator.

\medskip
{\bf Acknowledgments}. Both authors were supported by the Israel
Science Foundation grant ISF-248/08. We are grateful to Ofir Harari for
his assistance in running simulation examples and Saharon Rosset for
fruitful discussions.

\section*{Appendix}
Throughout the proofs we use $C$ to denote a generic positive constant, not
necessarily the same each time it is used, even within a single equation.
Similarly, $C(\gamma)$ is a generic positive constant depending on $\gamma$.

\subsection*{Proof of Theorem \ref{th:upper}}
As we have mentioned in Section \ref{sec:bayes}, the sparse group MAP estimator
can be viewed as a penalized likelihood estimator \fr{eq:est} with the
complexity penalties \fr{eq:penaltyj} and \fr{eq:penalty0}. We first
re-write it in a somewhat different form that will allow us then to apply
the general results of Birg\'e \& Massart (2001) for complexity penalized
estimators.

Let $\by=(y_{11},...,y_{n1},...,y_{1m},...,y_{nm})'$ be an
amalgamated $n \times m$ vector of data. Similarly,
$\bm=(\mu_{11},...,\mu_{n1},...,\mu_{1m},...,\mu_{nm})'$,
$\beps=(\epsilon_{11},...,\epsilon_{n1},...,\epsilon_{1m},...,\epsilon_{nm})'$
and the original model \fr{eq:model} can be re-written now as \be
\label{eq:model1} y_i=\mu_i+\epsilon_i,\;\;\;\epsilon_i
\stackrel{i.i.d.} \sim {\cal N}(0,\sigma_n^2),\;i=1,...,nm \ee
Define an indicator vector ${\bf d}$, where $d_i=\mathbb{I}\{\mu_i
\ne 0\},\; i=1,...,nm$. In terms of the model \fr{eq:model1},
$h_j=\sum_{i=n(j-1)+1}^{nj}d_i,\;j=1,...,m$ and $m_0=\#\{j: h_j >
0\}$. For a given ${\bf d}$, define $D_{\bf d}=\sum_{j=1}^m
h_j=\#\{i: d_i=1,\;i=1,...,nm\}$ and
$$
L_{\bf d}=\frac{1}{D_{\bf d}}\left(\sum_{j=1}^m \ln\left(\pi_j^{-1}(h_j)
{n \choose h_j} \right)+\ln\left(\pi_0^{-1}(m_0){m \choose m_0}\right)\right)
$$
for ${\bf d} \not \equiv {\bf 0}$ and $L_{\bf 0}=2\ln \pi_0^{-1}(0)$, where
we formally set $\pi_j(0)=1$.
Then, the sparse group MAP estimator
$\bhm=(\hat{\mu}_{11},...,\hat{\mu}_{n1},...,\hat{\mu}_{1m},...,\hat{\mu}_{nm})'$ is the penalized
likelihood estimator of $\bm$ with the complexity penalty
\begin{eqnarray*}
Pen({\bf d})& = & 2\sigma_n^2(1+1/\gamma)\left(\sum_{j=1}^m
\ln\left(\pi_j^{-1}(h_j) {n \choose h_j}
(1+\gamma)^{\frac{h_j}{2}}\right)+
\ln\left(\pi_0^{-1}(m_0){m \choose m_0}\right)\right) \\
& = & \sigma_n^2(1+1/\gamma)D_{\bf d} \left(2L_{\bf
d}+\ln(1+\gamma)\right)
\end{eqnarray*}
for ${\bf d} \not \equiv {\bf 0}$ and $Pen({\bf
0})=\sigma_n^2(1+1/\gamma)L_{\bf 0}$.

One can verify that
$$
\sum_{{\bf d} \not \equiv {\bf 0}}e^{-D_{\bf d} L_{\bf d}}=\sum_{k=1}^m
\pi_0(k)=1-\pi_0(0)
$$
A straightforward calculus (see the proof of Theorem 1 of
Abramovich, Grinshtein \& Pensky, 2007 for more details) implies
also that for any ${\bf d}$ under Assumption (P),
$$
(1+1/\gamma)(2L_{\bf d}+\ln(1+\gamma)) \geq C(\gamma) (1+\sqrt{2 L_{\bf d}})^2,
$$
where $C(\gamma)>1$.
One can then apply Theorem 2 of Birg\'e \& Massart (2001) to get
\begin{eqnarray}
\sum_{j=1}^m E||\bhm_j-\bm_j||^2_2 & \leq & c_1(\gamma) \min_{\J \subseteq \{1,...,m\}} \left\{
\sum_{j \in \Jc}
\min_{1 \leq h_j \leq n}
\left(\sum_{i=h_j+1}^n \mu^2_{(i)j}+Pen_j(h_j)\right) \right. \nonumber \\
& + & \left. \sum_{j \in \J}\sum_{i=1}^n \mu^2_{ij}+
Pen_0(m_0)\right\} +c_2(\gamma) \sigma_n^2 (1-\pi_0(0))
\label{eq:upper1}
\end{eqnarray}
\newline $\Box$

\subsection*{Proof of Theorem \ref{th:upperth}}
One can easily check from Table \ref{tab:rates} that for 
$\eta_{jn} \geq n^{-1/\min(p_j,2)}\sqrt{\ln n}$ for $p_j>0$, 
the last term $c_2(\gamma)\sigma_n^2(1-\pi_0(0))$ in the RHS of \fr{eq:upper}
is of order
$O(\sigma_n^2)=o(R(\Theta_j[\eta_{jn}]))$ for all nonzero $\bm_j$ and all
$p_j \geq 0$. 

Let $\J^{c*}$ be the true (unknown) subset of nonzero $\bm$'s and
$m^*_0=|\J^{c*}|$.

\medskip
\noindent
I. $m_0^* \leq \lfloor m/e \rfloor$.
\newline
Apply Theorem \ref{th:upper} for $\J=\J^*$:
\begin{eqnarray*}
\sum_{j=1}^m E||\bhm_j-\bm_j||^2_2 & \leq &
c_1(\gamma)\left\{\sum_{j \in {\cal J}_0^{c*}} \min_{1 \leq h_j \leq
n} \left(\sum_{i=h_j+1}^n \mu_{(i)j}^2
+2\sigma_n^2(1+1/\gamma)\ln\left(\pi^{-1}_j(h_j){n \choose
h_j}(1+\gamma)^{\frac{h_j}{2}}\right)\right) \right. \\
& + & \left. 2\sigma_n^2(1+1/\gamma)\ln\left(\pi_0^{-1}(m_0){m
\choose m_0}\right)\right\} +c_2(\gamma)\sigma_n^2(1-\pi_0(0))
\end{eqnarray*}
Since for $m_0=1,...,\lfloor m/e \rfloor$, ${m \choose m_0} \leq
(m/m_0)^{2m_0}$ (see Lemma A1 of Abramovich {\em et. al}, 2010), the
required conditions on $\pi_0(\cdot)$ ensure that
$$
2\sigma_n^2(1+1/\gamma)\ln\left(\pi_0^{-1}(m_0){m \choose
m_0}\right) \leq C(\gamma) \sigma_n^2 m_0\ln(m/m_0)
$$
To complete the proof for this case we consider now separately 
\be
\label{eq:minhj} \min_{1 \leq h_j \leq n} \left(\sum_{i=h_j+1}^n
\mu_{(i)j}^2 +2\sigma_n^2(1+1/\gamma)\ln\left(\pi^{-1}_j(h_j){n
\choose h_j}(1+\gamma)^{\frac{h_j}{2}}\right)\right) 
\ee 
for each $j \in {\cal J}_0^{c*}$ and show that it is
$O(R(\Theta_j[\eta_{jn}]))$ (see Table \ref{tab:rates}). We
distinguish between several cases, where the proofs for strong
$l_p$-balls will follow immediately from the proofs for the
corresponding weak $l_p$-balls due to the embedding properties
mentioned in Section \ref{sec:main}.

\medskip
\noindent {\em Case 1}: $\bm_j \in \Theta_j[\eta_{jn}],\;\eta_{jn} >
e^{-c(\gamma)}$ for $p_j=0$ and $\eta^{p_j}_{jn} > e^{-c(\gamma)}$
for $p_j>0$. Taking $h^*_j=n$, under the condition on $\pi_j(n)$
implies that \fr{eq:minhj} is $O(\sigma_n^2n)=O(R(\Theta_j[\eta_{jn}]))$.

\medskip
\noindent {\em Case 2}: $\bm_j \in l_0[\eta_{jn}],\; \eta_{jn} \leq
e^{-c(\gamma)}$. Note that since $\bm_j \not \equiv {\bf 0}$,
$\eta_{jn} \geq n^{-1}$. Choose $h^*_j=n\eta_{jn}$ and repeat the
arguments of the proof of Theorem 3 of Abramovich, Grinshtein \&
Pensky (2007) using a slightly more general Lemma A1 of Abramovich
{\em et. al} (2010) for approximating the binomial coefficient in
\fr{eq:minhj} instead of their original Lemma A.1.


\medskip
\noindent {\em Case 3}: $\bm_j \in m_{p_j}[\eta_{jn}],\; 0 < p_j < 2,\;
n^{-1}(\ln n)^{p_j/2} \leq \eta^{p_j}_{jn} \leq e^{-c(\gamma)}$. Take
$1 \leq h^*_j=n\eta^{p_j}_{jn}(\ln \eta^{-p_j}_{jn})^{-p_j/2} \leq ne^{-c(\gamma)}$ and follow
the proof of Theorem 4 of Abramovich, Grinshtein \& Pensky (2007)
with a more general version of Lemma A1 (see Case 2).

\medskip
\noindent {\em Case 4}: $\bm_j \in m_{p_j}[\eta_{jn}],\; p_j \geq 2,\;
n^{-p_j/2}(\ln n)^{p_j/2} \leq \eta^{p_j}_{jn} \leq e^{-c(\gamma)}$.
Take
$h^*_j=1$. Then, for $p_j>2$
$$
\sum_{i=h^*_j+1}^n \mu^2_{(i)j} <  \sigma_n^2 n^{2/p_j} \eta^2_{jn}
\int_1^nx^{-2/p_j}dx < \frac{p_j}{p_j-2} \sigma_n^2 n^{2/p_j}
\eta_{jn}^2 n^{1-2/p_j}=O(\sigma_n^2 n \eta^2_{jn})
$$
and, similarly, for $p_j=2$
$$
\sum_{i=h^*_j+1}^n \mu^2_{(i)j} <  \sigma_n^2 n \eta^2_{jn}
\int_1^nx^{-1}dx = \sigma_n^2 n \eta^2_{jn} \ln n
$$
On the other hand, under the conditions on $\pi_j(\cdot)$, $\pi_j(1)
\geq n^{-c_1}$ that yields
$$
2\sigma_n^2(1+1/\gamma)\ln\left(\pi^{-1}_j(1)n
\sqrt{1+\gamma}\right)=O(\sigma_n^2 \ln n) = O(\sigma_n^2 n
\eta^2_{jn})
$$
for $\eta_{jn} \geq \sqrt{n^{-1}\ln n}$.

\medskip
\noindent
II. $\lfloor m/e\rfloor < m^*_0 \leq m$.
\newline
Apply Theorem
\ref{th:upper} for $\Jc = \{1,...,m\}$ (or, equivalently, $\J =
\emptyset$) and $h_j=1$ for $j \in \J^*$~:
\begin{eqnarray}
\sum_{j=1}^m E||\bhm_j-\bm_j||^2_2 & \leq & c_1(\gamma)\left\{\sum_{j \in
{\cal J}_0^{c*}}\min_{1 \leq h_j \leq n}
\left(\sum_{i=h_j+1}^n \mu^2_{(i)j}+Pen_j(h_j)\right) +
\sum_{j \in {\cal J}_0^*}Pen_j(1)+Pen_0(m) \right\} \nonumber \\
& + & c_2(\gamma)\sigma_n^2(1-\pi_0(0)),
\label{eq:m0m}
\end{eqnarray}
where the conditions on $\pi_j(1)$ and
$\pi_0(m)$ imply $\sum_{j \in {\cal J}_0^*}Pen_j(1)= O(\sigma_n^2 m
\ln n)$ and $Pen_0(m)=O(\sigma_n^2 m)$. From Table \ref{tab:rates} one can
verify that for all dense and sparse cases, 
$\sigma^2_n \ln n=O(R(\Theta_j [\eta_{jn}]),\;j \in {\cal J}_0^{c*}$
and, therefore,
the first term $\sum_{j \in {\cal J}_0^{c*}}$ in the RHS of \fr{eq:m0m} 
is dominating for $m^*_0 \sim m$.
\newline $\Box$

\subsection*{Proof of Theorems \ref{th:lower1}-\ref{th:lower2}}
The ideas of the proofs of both theorems on the minimax lower bounds are similar
and can be combined.

Note first that any estimator cannot perform better than an oracle
who knows the true $\J$. In this (ideal) case one would obviously set
$\bhm_j \equiv {\bf 0}$ for all $j \in \J$ with zero risk and, therefore,
due to the additivity of the risk function,
$$
\inf_{\btm_1,...,\btm_m}
\sup_{\bm_j \in \Theta_j[\eta_{jn}],j \in \Jc}\sum_{j=1}^m E||\btm_j-\bm_j||^2_2
\geq C \sum_{j \in \Jc} R(\Theta_j[\eta_{jn}])
$$
for any $\Theta_{jn}[\eta_{jn}]$ (see, e.g., Johnstone, 2011, Proposition 4.14).

Furthermore, following Case II in the proof of  Theorem \ref{th:upperth}, 
$\sum_{j \in \Jc} R(\Theta_j[\eta_{jn}])$ dominates over 
$\sigma_n^2 m_0\ln(m/m_0)$ in \fr{eq:lower1} and \fr{eq:lower2} for $m_0>m/2$.
To complete the proof we need to show, therefore, that for $m_0 \leq m/2$,
the minimal unavoidable
price for not being an oracle for selecting nonzero $\bm_j$'s is of
order $\sigma_n^2 m_0\ln(m/m_0)$. 

The main idea of the proof is
to find a subset ${\cal M}_{m_0}$ of $n \times m$ vectors
$\bm=(\mu_{11},...,\mu_{n1},...,\mu_{1m},...,\mu_{nm})'$ with 
$m_0$ nonzero $\bm_j=(\mu_{1j},...,\mu_{nj})' \in
\Theta_j[\eta_{jn}]$ such that for any pair $\bm^1,\bm^2 \in {\cal
M}_{m_0}$ and some $C>0$, $||\bm^1-\bm^2||^2_2 \geq
C\sigma_n^2 m_0\ln(m/m_0)$, while the Kullback-Leibler
divergence $K(\mathbb{P}_{\bm^1},\mathbb{P}_{\bm^2})=
||\bm^1-\bm^2||^2_2/(2\sigma_n^2) \leq (1/16)\ln{\rm card}({\cal
M}_{m_0})$. The result will then follow immediately from Lemma A.1
of Bunea, Tsybakov \& Wegkamp (2007). 

Define the subset ${\cal \tilde{D}}_{m_0}$ of all $m$-dimensional indicator
vectors with $m_0$ entries of ones, that is
${\cal \tilde{D}}_{m_0}=\{\bd: \bd \in \{0,1\}^m,\; ||\bd||_0=m_0\}$. 
By Lemma A.3 of Rigollet \& Tsybakov (2011), for $m_0 \leq m/2$  there exists a subset
${\cal D}_{m_0} \subset {\cal \tilde{D}}_{m_0}$ such that for some constant
$\tilde{c}>0$, $\ln {\rm card}({\cal D}_{m_0}) \geq \tilde{c}m_0\ln(m/m_0)$,
and for any pair $\bd_1,\bd_2 \in {\cal D}_{m_0}$, the Hamming distance
$\rho(\bd_1,\bd_2)=\sum_{j=1}^m \mathbb{I}\{\bd_{1j} \neq \bd_{2j}\} \geq \tilde{c}
m_0$.

To any indicator vector $\bd \in {\cal D}_{m_0}$ assign the
corresponding mean vector $\bm \in {\cal M}_{m_0}$ as follows. Let
$\tilde{C}^2=(1/16)\sigma_n^2\tilde{c} \ln(m/m_0)$. Define
$\bm_j=(\tilde{C},0,...,0)'\mathbb{I}\{d_j=1\}$ for $0 \leq p_j < 2$
and
$\bm_j=(\tilde{C}n^{-1/2},\tilde{C}n^{-1/2},...,\tilde{C}n^{-1/2})'\mathbb{I}\{d_j=1\}
$ for $p_j \geq 2,\;j=1,...,m$. Hence, ${\rm card}({\cal M}_{m_0})={\rm
card}({\cal D}_{m_0})$. Obviously, the resulting $\bm_j \in
l_0[\eta_{jn}]$ and a straightforward calculus shows that under the
additional constraint on $\eta_{jn}$ and $m_0$ in Theorem
\ref{th:lower2}, $\bm_j \in l_{p_j}[\eta_{jn}]$.

For any $\bm^1, \bm^2 \in {\cal M}_{m_0}$ and the corresponding
$\bd_1, \bd_2 \in {\cal D}_{m_0}$, we then have
$$
||\bm^1-\bm^2||^2_2 = \tilde{C}^2 \sum_{j=1}^m
\mathbb{I}\{\bd_{1j}\neq\bd_{2j}\} \geq \tilde{C}^2 \; \tilde{c}\;
m_0 = (1/16)\sigma_n^2 \tilde{c}^2 m_0\ln(m/m_0)
$$
and
$$
K(\mathbb{P}_{\bm^1},\mathbb{P}_{\bm^2}) =
\frac{\tilde{C}^2}{2\sigma_n^2} \sum_{j=1}^m
\mathbb{I}\{\bd_{1j}\neq\bd_{2j}\} \leq \frac{\tilde{C}^2
m_0}{\sigma_n^2} \leq (1/16)\ln{\rm card}({\cal M}_{m_0})
$$
$\Box$

\end{document}